\documentclass[english,a4paper,oneside,11pt]{article}
\usepackage{babel}
\usepackage[dvips]{graphicx}
\newtheorem{theorem}{Theorem}[section]
\newtheorem{proposition}[theorem]{Proposition}

\newcommand{\titulo}[1]{\mbox{} \\ \noindent \textit{\textbf{\large #1}}\\}
\newcommand{\autor}[1]{\noindent \textbf{#1}}
\newcommand{\afil}[1]{{\small \noindent \textit{#1}}}
\renewcommand{\abstract}[1]{{\small \noindent \textbf{Abstract:} #1\\}}

\newcommand{\keywords}[1]{{\small \noindent \textbf{Keywords:} #1\\}}

\begin{document}
 \thispagestyle{empty}

\titulo{Nearly Gaussian random variables and application to meteorology}
\autor{Rui Gon\c{c}alves }

\afil{LIAAD - INESC TEC and Faculty of Engineering of the University of Porto, R. Dr. Roberto Frias, 4200-465 Porto, rjasg@fe.up.pt. Telephone/fax: 351225081921, Cellular tel.: 351963487610}

\abstract{
We consider a nearly Gaussian random variable $X$ (see \cite{Lefebvre}) that, after a power transformation, the variable $X^c$
where $c=\{(2k+1)/(2j+1)\}$, $k,j=\{0,1,\dots\}$, is approximately Gaussian. This transformation is useful to model errors in temperature forecasts.}


\keywords{kurtosis, power transformation, Meteorology, nearly Gaussian.}



\section{Introduction}
\label{sec1}
Let $x_1, x_2,\dots, x_n$ be a random sample of the random variable $X$. We define:
$$
\bar{x}=\sum_{i=1}^n\frac{x_i}{n}
$$
so that,
$$
\hat{\mu}_k=\sum_{i=1}^n\frac{(x_i-\bar{x})^k}{n-1}
$$
for $k=1,2,...$. The sample standard deviation $s_x=\sqrt{\hat{\mu_2}}$. The sample skewness coefficient is
$$
b_1=\frac{\hat{\mu}_3}{s_x^3}
$$
and the sample kurtosis is,
$$
b_1=\frac{\hat{\mu}_4}{s_x^4}
$$
A random variable having sample skewness and kurtosis near 0 and 3, respectively, is called quasi or nearly Gaussian, see Lefebvre \cite{Lefebvre}. If a random variable is the sum of independent and identically distributed (i.i.d.) random variables with finite mean and variance then by the central limit theorem it is nearly Gaussian.
To certify how close a random variable is to normality, one should use several goodness-of-fit tests, e. g., Qui-Square, Lilliefors and Shapiro-Wilk, etc. In case of rejection one can fit another symmetrical probability density function (pdf) or, instead, one can consider a simple power transformation of data, see Lefebvre \cite{Lefebvre}, and for examples of the use of power transformations in different contexts see Gon\c calves et al. \cite{GonPinSto1}, \cite{GonFerPinSto2} and \cite{GonFerPinSto3} .
Let $X$ be a nearly gaussian random variable, then we define,
\begin{equation}
Y=X^c
\label{ref1}
\end{equation}
where $c:=(2k+1)/(2j+1)$, $k,j\in\{0,1,2,...\}$ and $c$ close to 1. In that case, $Y$ is a power transformation of $X$. This transformation is bijective so that all data, negative and positive, can be equally transformed. Moreover, the power transformation is related to the well known Box-Cox transformation (see Box and Cox \cite{BoxCox}). However, an important difference between the two is that the power transformation does not affect symmetry. That is not the case of the Box-Cox, because generically the Box-Cox transformation destroys, if present, the symmetry of the data.
If the Gaussian distribution is acceptable for the transformed data ($Y$) then we can say that
$$
X=Y^{1/c}
$$
is approximately a power transformation of a Gaussian distribution.\\
In section \ref{sec2} we compute the kurtosis coefficient of the variable Y when X is exactly gaussian N(0,1).
We also show that based in the sample kurtosis coefficient and in a kurtosis table relating the power exponent with the theoretical kurtosis we can estimate what $c$ is appropriate to transform the data to normality.
In section \ref{sec3} we apply the method described in section \ref{sec2} to the data consisting of one day ahead forecast errors in daily maximum and minimum temperatures for the year 2011\footnote{The data was collected by Instituto Portugu\^es do Mar e da Atmosfera (IPMA)}.  Firstly, we consider the errors in maximum temperatures. Using the Lilliefors and the Shapiro-Wilk test we conclude that normality is rejected. Then using the sample kurtosis and the table \ref{tab1} we found the exponents $9/11$ and $7/9$ to be appropriated to transform the original data in a sample with Gaussian distribution. This means that the original data can be well fitted by a Gaussian distribution raised to the power $11/9$ or $9/7$. In order to compare power normal with other symmetric distributions we perform a Qui-square goodness-of-fit test for the power normal and for the Laplace and Pearson IV distributions. The pvalues observed are greater than the usual significance levels but the pvalue of the power normal is significantly greater than that of the Laplace and Pearson IV distribution. In the case of the one day ahead forecast of minimum temperature errors we found that normality is not reject.

\section{Kurtosis coefficient of the power transformation of a Gaussian random variable}
\label{sec2}

Let $Z$ be a Gaussian random variable with parameters $\mu=0$ and variance $\sigma^2$, then for $c>0$,
$$
E(Z^c)=\int_{-\infty}^{\infty}\frac{z^c}{\sqrt{2\pi}\sigma}=\frac{2^{c/2}}{2\sqrt{\pi}}[1+(-1)^c]\Gamma\left(\frac{c}{2}+\frac{1}{2}\right)
$$
where $\Gamma$ is the gamma function. Hence, the following proposition may be stated,
\begin{proposition}
The kurtosis coefficient of the random variable defined in (\ref{ref1}) when $X\sim N(0,\sigma^2)$ is given by,
$$
\beta_2(c)=\frac{\int_{-\infty}^{\infty}\frac{1}{\sqrt{2\pi}}x^{4c}e^{-x^2/2}dx}{\left(\int_{-\infty}^{\infty}\frac{1}{\sqrt{2\pi}}x^{2c}e^{-x^2/2}dx\right)^2}=\sqrt{\pi}\frac{\Gamma(2c+\frac{1}{2})}{\Gamma^2(c+\frac{1}{2})}
$$
\end{proposition}
\begin{table}[htb]
\centering
\begin{tabular}{|c|c|c|c|}
 \hline
  $c$      & $\beta_2(c)$ & $c$ & $\beta_2(c)$\\
   \hline
  3/5 & 1.779 & 15/11 & 4.894\\
  5/7 & 2.06 & 7/5 & 5.142\\
  7/9 & 2.237 & 17/15 & 3.58\\
  9/11 & 2.358 & 15/13 & 3.08\\
  11/13& 2.447 & 11/15 & 2.11\\
  13/15& 2.514& 19/21& 2.643\\
  15/17& 2.566& 23/25& 2.697\\
  17/19& 2.6  & 29/31& 2.753\\
  13/11& 3.828& 31/33& 2.767\\
  23/19& 3.979& 21/23& 2.672\\
  11/9& 4.042& 35/37& 2.79\\
  9/7 & 4.404& 37/39& 2.8\\
 \hline
\end{tabular}
\caption{Kurtosis of the random variable $Y=X^c$ where $X\sim N(0,1)$ for a few values of $c$.}
\label{tab1}
\end{table}
We present the table \ref{tab1} with a few values of the quantity $\beta_2(c)$ for some values of $c$. This table will be used to select an exponent $c$, close to 1, so that, raising the values of the nearly Gaussian sample to $1/c$ one obtains an approximately Gaussian sample.
Note that the identification of the exponent $c$ requires that the variable should have mean and skewness close to zero. To apply this transformation to non-symmetrical data, a Box-Cox transformation to symmetry should be applied first, see Hinkley \cite{Hinkley}.
Given a sample $x_1,\dots,x_n$ we center the data defining
$$
z_i=x_i-\bar{x}
$$
so that the mean of the $z_i$ is 0. Now assuming that the skewness coefficient is close to zero and the kurtosis is not close to 3 to be Gaussian then, by selecting an appropriate $c$ from table \ref{tab1}, it is possible to find the transformation $W=Z^{1/c}$ that is likely to transform the data into a sample that is approximately Gaussian.
\begin{proposition}
If $Y=X^c$ is gaussian, $Y\sim N(\mu,\sigma^2)$, then pdf of the power transformation of a Gaussian variable, $X=Y^{1/c}$ is,
$$
f_X(x)=\frac{1}{\sqrt{2\pi}c}c|x^{c-1}|\exp{\left(-\frac{x^c-\mu}{2\sigma^2}\right)}
$$
\end{proposition}
In \cite{Lefebvre} the author shows that the power transformation of a symmetric non Gaussian pdf such as the Laplacian can be very similar to that of a Gaussian. He supported the claim by conducting a simulation study where is shown that applying an adequate exponent to Laplacian samples, in 82 out of 100 samples, tested for normality, the null hypothesis was accepted.

\section{Application to temperature forecasts}
\label{sec3}
Our goal is to find statistical models for the forecasting errors of minimum and maximum temperatures. In \cite{Lefebvre} the author considered this problem and found out that the pdf of a power of a Gaussian random variable fitted with success the temperature error forecasts.
Here, we consider in first place the one day ahead maximum temperature forecasts during the year 2011. The size of the data set is 347 (there are a few missing data). We define $X$ as the forecasting error $T_F-T_O$ where $T_F$ is the forecast and $T_O$ the observed temperatures. The observed statistics for the error in the maximum temperature are:
\begin{center}
$\bar{x} =-0.0458$; $s_x$= 1.616
$b_1=0.035$; $b_2=3.567$
\end{center}
where $\bar{x}$, $s_x$, $b_1$ and $b_2$\footnote{some authors define the kurtosis as $b_2^*=b_2-3$ so that $b_2^*=0$ in the case of the gaussian.} are the sample mean, sample standard deviation, sample skewness coefficient and sample kurtosis coefficient, respectively.
\begin{figure}
\begin{center}
\includegraphics[width=10cm]{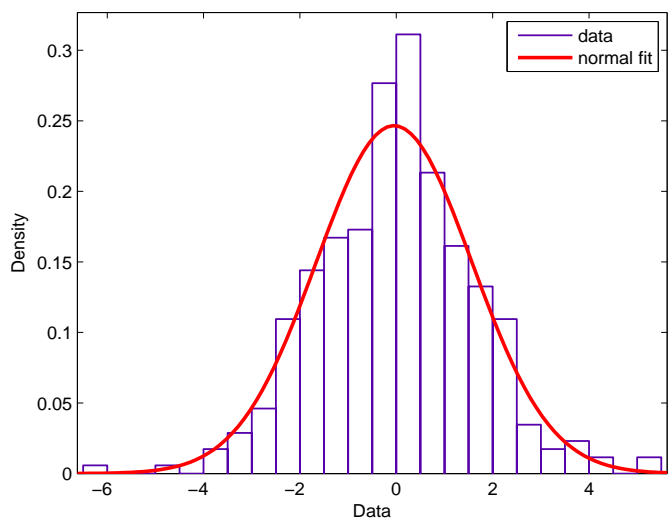}
\caption{Histogram and Gaussian fit to the errors in maximum temperature forecasts.}
\label{fig1}
\end{center}
\end{figure}
\noindent
In figure \ref{fig1} we present the histogram of the data with the Gaussian fit superimposed. As in Lefebvre \cite{Lefebvre} we found more observations in the center when compared to the Gaussian density. Before applying any statistical test we must make sure that there is none or very little temporal correlation among the data. A way of measuring temporal correlation is by computing the sample autocorrelation function (ACF). If the data is comparable to a white noise sequence then its ACF should be in the $95\%$ confidence bounds. We found 3 values in 20 outside the bounds which means that less than $95\%$ is inside the bounds. This implies that there is temporal correlation in the data. Hence, we decided to take a subset of the data. We kept only one in each pair of forecasting errors. This procedure reduces the data size to an half (173) which is still large enough. For the reduced data set the values of the ACF are inside the $95\%$ confidence bounds which means that the data is comparable to white noise. Using the software SPSS we performed the Lilliefors and Shapiro-Wilk tests and we obtained pvalues of 0.023 and 0.1 respectively. These quantities are less that the usual significance levels used in statistical tests. So, we conclude that the Gaussian distribution is not a good model for $X$. Next, we tried classic models for the forecast errors such as the laplacian and the Student's $T$ and in both cases it turned out that the models were not acceptable.
To look for an appropriate power transformation we must compute some statistics of the data first. The mean is close to zero so there is no need to center the data. The kurtosis of the reduced data set is 4.0583. Based on the kurtosis in table \ref{tab1}, we tried the transformation,
$$
 w_k=x_k^{9/11}
$$
We found that the gaussian distribution is acceptable as a model to the data. Applying the Lilliefors test (with Matlab) the p-value increased from 0.023 to at least 0.2 and the p-value of the Shapiro-Wilk test is 0.439. Because the values of the sample kurtosis are not exactly equal to those of the table we also tried another exponent close to 9/11. In fact, trying,
$$
 w_k=x_k^{7/9}
$$
and applying the Lilliefors and Shapiro-Wilks tests (with SPSS) we obtain the p-values 0.2 and 0.439 which are equal to those obtained for the former exponent.
Hence, we can use a Gaussian distribution raised to the power $11/9$ to model the data. Next, raising the original data to the power $9/11$ we obtain an approximately Gaussian distribution with mean -0.027 and standard deviation 1.311,
$$
X^{(9/11)} \approx N(-0.027,1.311^2)\nonumber.
$$
The exponents that we have found are different from the ones reported in Lefebvre \cite{Lefebvre}. In the same paper, the author fitted classical distributions such as the Laplace distribution but found that none of them was acceptable. Because our data is close to symmetry we considered the power normal, the laplace distribution and the Pearson IV distributions. We performed a Chi-Square goodness-of-fit test. The results are shown in table \ref{tab2}.
Firstly, we consider the Laplace distribution,
$$
f(x|\mu,b) = \frac{1}{2b} \exp \left( -\frac{|x-\mu|}{b} \right) \,\!
$$
The estimate for $\mu$ is zero and the maximum likelihood estimate for $b=\sum_{i=1}^n|X_i-\hat{\mu}|$.
Secondly, we considered the Pearson type IV distribution, its pdf is,
$$
f_X(x)= k\left[1+\left(\frac{x-\lambda}{a}\right)^{-m}\right]\exp\left\{-\nu\arctan{\left(\frac{x-\lambda}{a}\right)}\right\}\nonumber
$$
for $x\in \Re$ and the parameters $m$,  $m>1/2$, $\nu$, $a$ and $\lambda$ are real constants. The normalizing constant $k$ is given by,
$$
k=\frac{2^{2m-2}|\Gamma[m+(i\nu/2)]|^2}{\pi a\Gamma(2m-1)}\nonumber
$$
Using the estimators given by the method of moments (see \cite{StuartOrd}) we obtained,
$$
\hat{m}= 5.5666, \hat{\nu}=-1.6843, \hat{a}=4.2831, \hat{\lambda}= -0.8211
$$
The results of the chi-square tests are presented in table 4.
\begin{table}[htb]
\centering
\begin{tabular}{|l|c|c|c|c|}
\hline
interval       &$n_j$	&	$X$   	    &	Lap	        &	Pearson IV	\\
\hline
($-\infty$, -2.5)&	8	&	10.6870	&	9.4782	&	4.4115	\\
(-2.5, -2)	     &	11	&	7.01236	&	5.2715	&	7.9234	\\
(-2, -1.5)	     &	12	&	10.4451 &	8.20341	&	12.8885	\\
(-1.5, -1,0)     &	11	&	14.8661	&	12.7660	&	18.4764	\\
(-1, -0.5)	     &	16	&	20.4931 &	19.8660	&	22.9571	\\
(-0.5, 0)	     &	31	&	32.2387	&	30.9149	&	24.566	\\
(0,0.5)	         &	31	&	25.1018	&	30.9149	&	22.7495	\\
(0.5, 1)	     & 	21	&	17.5664	&	19.8659	&	18.4937	\\
(1,1.5)	         &	10	&	12.4670 &	12.7659	&	13.4767	\\
(1.5,2)	         &	9	&	8.50495 &	8.20342	&	8.996	\\
(2, 2.5)	     &	7	&	5.54012	&	5.27154	&	5.6398	\\
(2.5, $\infty)$  &	6	&	7.92306	&	9.47825	&	3.3908	\\
$d^2$	         &		&	8.63832	&	14.0892	&	13.1724	\\
p-value	         &		&	0.47131	&	0.11919	&	0.15496	\\
\hline
\end{tabular}
\caption{Chi-Square goodness of fit test to determine  whether a Gaussian $N(-0.027, 1.311^2)$ distribution raised to the power $11/9$ is a good model for the raw data. The four columns give the chosen subintervals, the number $n_j$ of observations in each subinterval and the expected $e_j$ number of observations for each subinterval and for the Gauss, Laplace and Pearson IV distributions in this order.}
\label{tab2}
\end{table}
We see that the p-value of the observed statistic for the power transformation of the normal is considerably better than the others. We also checked by using a Lilliefors test, whether a Gaussian pdf fits the raw data raised to the power 9/11.
Next, we turn our attention to the forecasts errors of minimum temperatures. Like in the case of the maximum temperature when computing the ACF there are 2 values outside the $95\%$ confidence bounds  therefore we decided to eliminate one value in each two.
The observed statistics for the error in the minimum temperature for the reduced data set are:
\begin{center}
$\bar{x} =0.099$; $s_x$= 1.52
$b_1=-0.04$; $b_2=2.67$
\end{center}
Applying the Lilliefors and the Shapiro-Wilk tests (again with SPSS) we found a pvalue of at least 0.2 for the Lilliefors test and a pvalue of 0.439 for the Shapiro-Wilk test. This means that the data is already approximately Gaussian and there is no need for transformations.

\section{Concluding remarks}
\label{sec4}
In this article, we show that a using an appropriate exponent of the form $(2k+1)/(2j+1)$, ${k,j}=0,1,...$ the power transformation of a nearly Gaussian random variable can be Gaussian. The transformation is bijective so it may be used in both positive and negative data. We applied the method described in section \ref{sec2} to data consisting of the one day ahead forecast errors in daily maximum and minimum temperatures. In the case of errors in maximum temperatures we used both Lilliefors and Shapiro-Wilk tests and we concluded that normality was rejected. Afterwards, using the sample kurtosis and the table \ref{tab1} we found the appropriate exponents $9/11$ and $7/9$ that transform the original data in a sample with Gaussian distribution. Fitting results of the power normal, Laplace and Pearson IV distributions were compared and the pvalue of the power normal was found to be significantly greater than the other two. Surprisingly, in the case of the one day ahead forecast of daily minimum temperature errors, normality was not rejected.




\bibliographystyle{elsarticle-num}



\end{document}